\documentclass[final,pdftex,reqno,a4paper,12pt]{myectaart}

\RequirePackage[OT1]{fontenc}
\RequirePackage{amsthm}
\RequirePackage[cmex10]{amsmath}
\RequirePackage{natbib}
\RequirePackage[colorlinks,citecolor=blue,urlcolor=blue]{hyperref}
\RequirePackage{hypernat}

\usepackage{a4wide,graphicx,float,amssymb} 

\startlocaldefs
\numberwithin{equation}{section}
\theoremstyle{plain}

\endlocaldefs

\newcommand{\distr}{\stackrel{d}{=}}

\newtheorem{theorem}{Theorem}

\newtheorem{propos}[theorem]{Proposition}

\RequirePackage{ae,fancyvrb}
\DefineVerbatimEnvironment{Sinput}{Verbatim}{fontshape=sl}
\DefineVerbatimEnvironment{Soutput}{Verbatim}{}
\DefineVerbatimEnvironment{Scode}{Verbatim}{fontshape=sl}

\setkeys{Gin}{width=0.9\textwidth}

\begin{document} 

\begin{frontmatter}
\title{\bf The generalized lognormal distribution \protect\\ and the Stieltjes moment problem %\protect\thanksref{T1}
}
\runtitle{Generalized lognormal distribution and Stieltjes moment problem}

%\thankstext{T1}{Version: \today.}

\begin{aug}
\author{\fnms{Christian} \snm{Kleiber}%\thanksref{t1}
\ead[label=e1]{christian.kleiber@unibas.ch}}

%\thankstext{t1}{Version: \today.}

\runauthor{C. Kleiber}

%\thankstext{t1}{Some comment}

\address{Version: \today. \ Correspondence to: Christian Kleiber, Faculty of Business and Economics, %Quantitative Methods Unit, 
Universit\"at Basel, Peter Merian-Weg 6, 4002 Basel, Switzerland.
\printead{e1}}
\end{aug}

\begin{abstract}
This paper studies a Stieltjes-type moment problem defined by the generalized lognormal distribution, a heavy-tailed distribution with applications in economics, finance and related fields. It arises as the distribution of the exponential of a random variable following a generalized error distribution, and hence figures prominently in the EGARCH model of asset price volatility. Compared to the classical lognormal distribution it has an additional shape parameter. It emerges that moment (in)determinacy depends on the value of this parameter: for some values, the distribution does not have finite moments of all orders, hence the moment problem is not of interest in these cases. For other values, the distribution has moments of all orders, yet it is moment-indeterminate. Finally, a limiting case is supported on a bounded interval, and hence determined by its moments. For those generalized lognormal distributions that are moment-indeterminate Stieltjes classes of moment-equivalent distributions are presented. 



\end{abstract}

%\begin{keyword}
%\kwd{Generalized error distribution} 
%\kwd{generalized lognormal distribution, lognormal distribution, moment problem, size distribution, %statistical distributions, 
%Stieltjes class, volatility model}
%\end{keyword}

\end{frontmatter}

\section{Introduction}

The moment problem asks, for a given distribution with distribution function (CDF) $F$ with finite moments $m_k (F) ~=~ \int ^\infty _{-\infty} x^k \ \hbox{d}F(x)$ of all orders $k = 1, 2, \dots, $ whether or not $F$ is uniquely determined by the sequence of these moments. If $F$ is uniquely determined by this sequence,  $F$ or a random variable $X$ following this distribution are called
moment-determinate (for brevity, M-det); otherwise $F$ or $X$ are called moment-indeterminate (M-indet). Cases where the support of the distribution $F$ is the positive half-axis $\mathbb{R}^+ = [0, \infty)$ are called Stieltjes moment problems, cases where the support is the real line are called Hamburger moment problems, and cases where the support is a bounded interval are called Hausdorff moment problems.

The probably most widely known example of an M-indeterminate distribution is the lognormal distribution, first described by \cite{mom:Stieltjes:1894} in a non-probabilistic setting and further developed by \cite{mom:Heyde:1963}. The lognormal distribution is a basic model for describing size phenomena in economics and related fields \citep[see, e.g.,][]{mom:Kleiber+Kotz:2003}, including distributions of personal income, actuarial losses, or city sizes. It also arises in mathematical finance in the fundamental geometric Brownian motion model of asset price dynamics. 
Given the central role of the lognormal distribution in Stieltjes-type moment problems it is, therefore, of special interest to explore closely related distributions with respect to M-indeterminacy. Recently, \cite{mom:Lin+Stoyanov:2009} studied a generalization of the lognormal distribution derived from a skewed generalization of the normal distribution, finding that it is M-indeterminate for every value of the skewness parameter. The present paper explores a family of generalized lognormal distributions derived from a more classical symmetric generalization of the normal distribution, which compared to the normal distribution has an additional shape parameter. Like the classical lognormal distribution, this generalized version has been employed in financial economics as well as in modeling size distributions.

It turns out that this family of distributions sheds new light on the classical lognormal moment problem, in that M-determinacy now depends on the value of the shape parameter. Specifically, the family incorporates heavy-tailed distributions for which not all integer moments exist, moderately heavy-tailed distributions for which all moments exist yet the distributions are M-indeterminate, and, as a limiting case, a distribution with bounded support that is, therefore, determined by its moments. It also emerges that the classical lognormal distribution does not constitute an extreme case within the family: in the setting considered here, there exist more as well as less heavy-tailed M-indet distributions than the lognormal.

The paper is organized as follows: Section 2 provides some background on the generalized lognormal distribution. Section 3 contains a characterization of moment (in)determinacy for the family of generalized lognormal distributions in terms of their shape parameter, while Section 4 describes Stieltjes classes pertaining to the indeterminate cases. Section 5 concludes.

\section{The generalized lognormal distribution}

Being one of the basic distributions in probability and statistics, the normal distribution has triggered a number of generalizations. One such generalization is defined by the density

\begin{equation}\label{generror} 
f(y) ~=~ \frac{1}{2 \ r^{1/r} \ 
\sigma \ \Gamma (1+ 1/r)} ~ \exp \left\{ - \frac{1}{r \ \sigma^r} |
y - \mu | ^r \right\} , \quad - \infty < y < \infty, 
\end{equation} 
which includes the normal as the special case where $r=2$. Here $\mu \in \mathbb{R}$ is a location parameter and $\sigma \in \mathbb{R}^+$ is a scale parameter. The new parameter $r \in \mathbb{R}^+$ is a shape parameter  measuring tail thickness, with lower values of $r$ indicating heavier tails. The parameter $r$ plays a crucial role below.

This distribution is fairly widely known; however, it is known under different names in different fields and it was (re)discovered several times in different contexts. Specifically,  since $r=2$ yields the normal distribution and $r=1$ the Laplace distribution, the distribution~(\ref{generror}) is known both as a generalized normal distribution, in particular in the Italian language literature \citep{mom:Lunetta:1963, mom:Vianelli:1963}, and as a generalized Laplace distribution. It is also known as the normal distribution of order $r$, again especially in the Italian literature \citep[e.g.,][]{mom:Vianelli:1983}, and as the generalized error distribution, notably in econometrics and finance \citep[e.g.,][]{mom:Nelson:1991}. A further name is exponential power distribution \citep{mom:Box+Tiao:1973}, the name under which this distribution is presumably best known in the statistical literature. To the best of the author's knowledge, the generalized form (\ref{generror}) was first proposed  in a Russian journal by \cite{mom:Subbotin:1923}, who sought an axiomatic basis for a generalized form of Gauss's ``law of error.'' Hence the name Subbotin distribution is also in use, notably in econophysics \citep[e.g.,][]{mom:Alfarano+Milakovic+Irle:2012}. A multivariate generalization of (\ref{generror}) is the Kotz-type distribution \citep{mom:Kotz:1975}.

In what follows we sometimes set $\mu = 0$, since in the context of moment problems no extra generality is gained by including this location parameter. There exist different parameterizations of (\ref{generror}), notably regarding the scale parameter, but for the purposes of this paper the relevant parameter is $r$, so this complication shall be ignored below.

The \emph{generalized lognormal distribution} \citep{mom:Vianelli:1982a, mom:Vianelli:1983}, or perhaps \emph{logarithmic generalized normal distribution}, is less widely known than the generalized normal distribution. In fact, most of the currently available works are written in Italian and published in Italian journals and collected volumes that are often not easily available outside of Italy. A more accessible source may be \citet[][Ch.~4.10]{mom:Kleiber+Kotz:2003}, who summarize many basic properties. The distribution is defined as the distribution of $X = \exp (Y)$, where $Y$ follows eq.~(\ref{generror}), leading to the density

\begin{equation}\label{genlnpdf} 
f(x) ~=~ \frac{1}{2 \ x \ r^{1/r} \ 
\sigma \ \Gamma (1+ 1/r)} ~ \exp \left\{ - \frac{1}{r \ \sigma^r} |
\ln x - \mu | ^r \right\} , \quad 0 < x < \infty. 
\end{equation}

\begin{figure}[ht!]
\begin{center}
\includegraphics{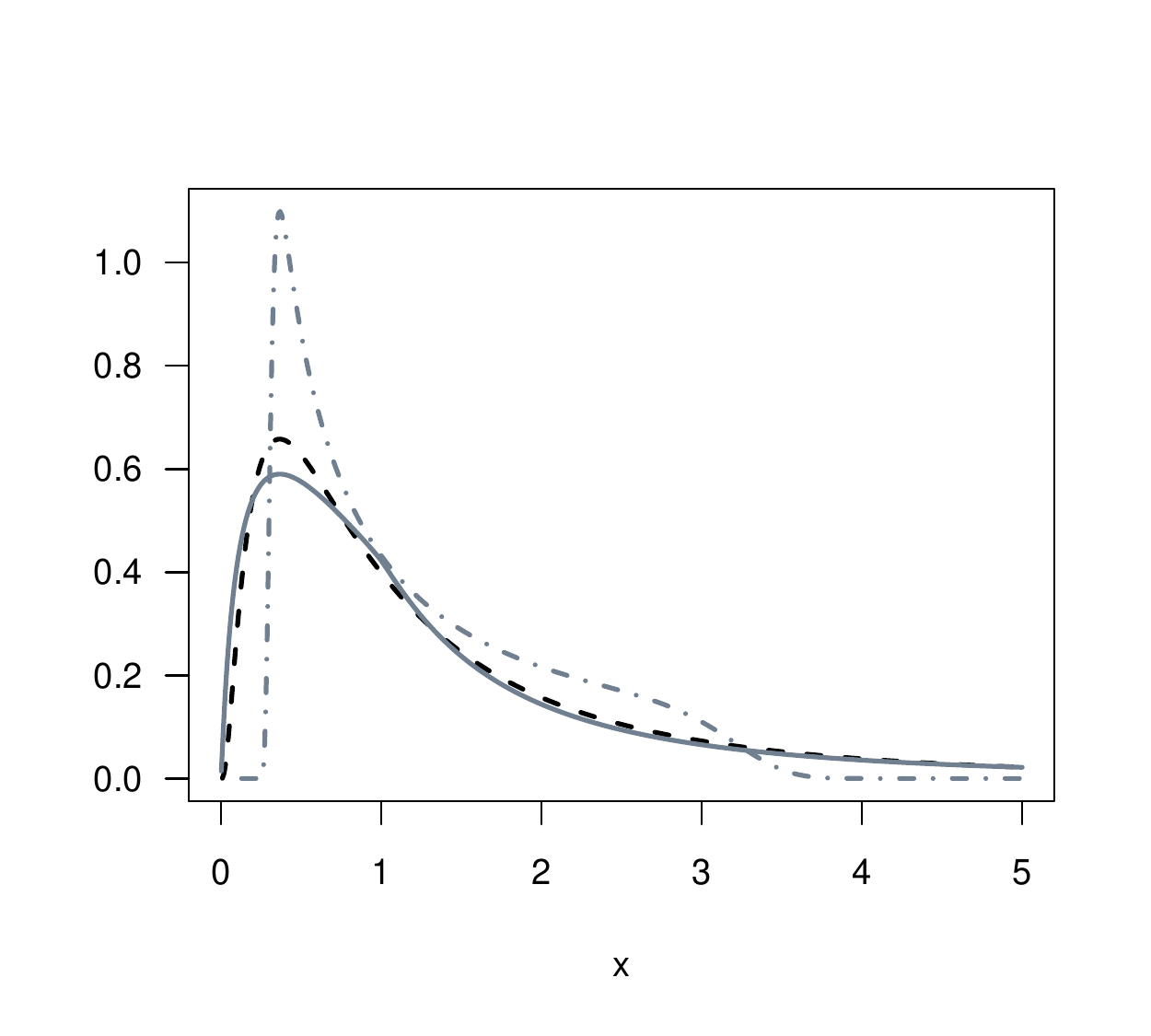}
\end{center}
\caption{\label{fig:fig1} Some generalized lognormal distributions (solid grey: $\mu = 0$, $\sigma = 1$, $r = 1.5$, %5, %10, 
dashed-dotted grey: $\mu = 0$, $\sigma = 1$, $r = 15$). The dashed black curve corresponds to the classical lognormal distribution ($r=2$, with $\mu = 0, \sigma = 1$). }
\end{figure}

If a random variable $X$ follows eq.~(\ref{genlnpdf}) this is denoted as $X \sim $ GLN($\mu, \sigma, r$). The distribution will sometimes be referred to as the generalized lognormal distribution of order $r$ if further emphasis is needed. The case where $r=2$ gives the classical lognormal distribution. In eq.~(\ref{genlnpdf}), $e^{\mu}$ is a scale parameter, while $\sigma$ and $r$ are both shape parameters. The effect of the new parameter $r$ is illustrated in Figure \ref{fig:fig1}. This Figure suggests that the density becomes more and more concentrated on a bounded interval with increasing $r$. Specifically, for $r=1.5$ the density is much like the classical lognormal density, but with slightly heavier tails, while for $r=15$ several points of inflection and a more rapid decrease in the tails emerge. The limiting case where $r \to \infty$ will also be explored below, see Theorem~\ref{thlimr}.

Like the classical lognormal distribution, the generalized lognormal distribution has been employed in economics and finance. As mentioned above, it has been used as a model for the size distribution of personal incomes. In an application to Italian income data, \cite{mom:Brunazzo+Pollastri:1986} estimate $r$ in the vicinity of 1.45, suggesting a model with even heavier tails than the classical lognormal distribution for their data. It will emerge below that their estimated model is not determined by its moments.

Perhaps more prominently, the distribution also arises in the widely used exponential GARCH (EGARCH) model of asset return dynamics~\citep{mom:Nelson:1991}, where it provides a more realistic specification of the innovation distribution in the volatility equation than the normal distribution. Recall that, in view of the exponential transformation employed in the EGARCH model, a widely used alternative to the normal distribution in GARCH modeling, the $t$ distribution, leads to tails that are too heavy, in the sense that the distribution corresponding to the exponentiated random variable has no moments of any order. In contrast, it will emerge below that the less extreme members of the generalized lognormal distribution possess moments of all orders, yet they are M-indeterminate. Specifically, all models estimated by \cite{mom:Nelson:1991}, with shape parameters $r$ in the vicinity of 1.56--1.57, are not determined by their moments. More recent work \citep[e.g.,][]{mom:Taylor:2005} confirms that $1 < r < 2$ is the empirically relevant range of the tail thickness parameter in this model. All of these objects are M-indeterminate.

\section{Generalized lognormal distributions and the moment problem}

How can one determine whether or not a given distribution with CDF $F$ is determined by the sequence of its moments? Although necessary and sufficient conditions are known \citep[see, e.g.,][]{mom:Shohat+Tamarkin:1950}, they are not very practical.  For M-determinacy, a sufficient condition is the existence of the moment generating function (MGF) $m_X (t) \ = \ \mathbb{E}[e^{tX}] \ = \ \int_0 ^\infty e^{tx} \ \hbox{d}F_X(x)$, $|t| < t_0$, for some $t_0 >0$. 

From the expression for the density (\ref{genlnpdf}) of the generalized lognormal distribution it is immediate that, for  any $r \in \mathbb{R}^+$, $\mathbb{E}[e^{t X}] = \infty$ for all $t > 0$; hence the MGF does not exist. It remains to explore the existence of the moments themselves. (Note that in view of $X > 0$ (a.s.) it is possible to consider moments of fractional order.) Without loss of generality, set $\mu = 0$ since $\exp(\mu)$ is a scale parameter. Substituting $z = \ln x$ yields, for some  $C > 0$,

\begin{equation}\label{exmom}
\mathbb{E}[X^k] ~=~ \int_0 ^\infty x^k f(x) \ \hbox{d}x ~=~ C \int_{-\infty} ^\infty \exp \{ k z -  |z|^r / (r \sigma^r) \} \ \hbox{d}z .
\end{equation}

This shows that convergence of the integral depends on the value of $r$: for $r > 1$ the integral is finite for all $k$, for $r = 1$ the condition $|k| < 1/\sigma$ is needed, while for $r < 1$ it does not converge for any $k \neq 0$. The following proposition collects these observations:

\begin{propos}\label{prop:exmom}
Suppose $X \sim GLN(\mu, \sigma, r)$.
\begin{enumerate}
\item The moment-generating function of $X$ does not exist for any $r \in (0, \infty)$.
\item The $k$th moment $\mathbb{E}[X^k]$ exists if and only if
\begin{itemize}
\item $k = 0$, if $r < 1$.
\item $|k| < 1/\sigma$, if $r = 1$.
\item $k \in (-\infty, \infty)$, if $r > 1$.
\end{itemize}
\end{enumerate}
\end{propos}

Apart from the integral representation (\ref{exmom}), it is also possible to obtain series expansions of the moments (when they exist). For $r >1$, they are of the form

\[
\mathbb{E}[X^k] \ = \   \frac{e^{k \mu}}{\Gamma\left( \frac 1r \right)} \ \sum_{i=0}^\infty \frac{(k \sigma )^{2i}}{(2i)!} \ r^{2i/r} \ \Gamma\left( \frac{2i+1}{r} \right), \quad k = 0,1,2, \dots,
\]
see \cite{mom:Brunazzo+Pollastri:1986} or \cite{mom:Nelson:1991}\footnote{It should be noted that these works employ different parameterizations of the distribution. Also, \cite{mom:Nelson:1991} obtains expectations of somewhat more general objects. Setting $\gamma = 0$, $p=0$ and $\theta=1$ in his Theorem A1.2 yields the required moments. The resulting expressions can be shown to coincide with those presented by \cite{mom:Brunazzo+Pollastri:1986}.}.

In view of Proposition \ref{prop:exmom} not all generalized lognormal distributions are of interest in the context of the moment problem. For $r=1$, only some moments exist, for $r < 1$ no moments exist. The cases where $r < 1$ therefore provide examples of distributions without any moments, integer or fractional. An earlier example was given by \cite{mom:Kleiber:2000}. For the remaining cases where $1 < r < \infty$ all the moments are finite yet the MGF does not exist. These are circumstances under which M-indeterminacy may arise. 

It remains to show that the distributions where $1 < r < \infty$ are indeed M-indet. For M-indeterminacy, a useful sufficient condition is the Krein condition \citep[e.g.][]{mom:Stoyanov:2000}. In a Stieltjes-type moment problem, it requires, for a density $f$ that is strictly positive for all $x \geq a > 0$, for some $a > 0$, that the normalized logarithmic integral of the density

\begin{equation}\label{krein}
K_S [f]  ~=~ \int_a ^\infty \frac{- \ln f(x^2)}{1 + x ^2} \ \hbox{d}x
\end{equation}
is finite. $K_S[f]$ is called the Krein integral of $f$. 

The following Theorem shows that generalized lognormal distributions of orders $1 < r < \infty$ are M-indeterminate:

\begin{theorem}
All generalized lognormal distributions $GLN(\mu, \sigma, r)$ of order $1 < r < \infty$ are M-indeterminate.
\end{theorem}

\noindent
{\sc Proof.} \ Setting without loss of generality $\mu = 0$ and $\sigma = 1$, the Krein integral (\ref{krein}) is, for $a > 0$ and $C_r > 0$ the normalizing constant,

\[
K_S[f] ~=~ \int_a ^\infty \frac {  -\ln C_r + 2 \ln x + \frac{1}{r} | 2 \ln x | ^r}{ 1 + x^2} \ \hbox{d}x .
\]

Since for large $x$ the integrand is eventually dominated by $x^{-1-\delta}$, for any $\delta \in (0,1)$, this integral is  finite for all $1 < r < \infty$, which gives the result. \hfill $\square$

\bigskip
Alternative proofs could employ results presented by \citet[][Remark~6.2]{mom:Gut:2002} or \citet[][p.~110]{mom:Pakes+Hung+Wu:2001}. 

\medskip
For $X_i \sim$ GLN($\mu_i, \sigma_i, r_i$), $i = 1,2$, with $r_i > 1$ and  densities $f_i$ it is easily seen that $\lim _{x \to \infty} f_1 (x)/ f_2(x) = \infty$ iff $r_1 < r_2$, hence the generalized lognormal distributions are, in a sense, ``more M-indeterminate'' for smaller $r$. (Indeed, in view of Proposition \ref{prop:exmom} for $r=1$ some moments no longer exist.) Specifically, the generalized lognormal distributions with $1 < r < 2$ are even more extreme than the classical lognormal distribution ($r=2$). Also, the cases where $2 < r < \infty$ are less extreme. It is also worth noting that although the tails of the generalized lognormal distribution become lighter and lighter with increasing $r$, the distribution is M-indet no matter how large $r$. It is, therefore, natural to ask what happens in the limit, i.e., for $r \to \infty$. The following Theorem addresses this case:

\begin{theorem}\label{thlimr}
For $r \to \infty$, the generalized lognormal distribution GLN($\mu,\sigma, r$) tends to a distribution supported on a bounded interval. Hence this limiting distribution is M-det.
\end{theorem}

\noindent
{\sc Proof.} \ It is convenient to analyze the limiting case for the distribution of $Y = \ln X$, i.e., the generalized normal distribution. Without loss of generality, set $\mu = 0$ and $\sigma = 1$. A random variable $Y$ following a generalized normal distribution admits the mixture representation \citep[][p.~175]{mom:Devroye:1986}

\begin{equation}\label{subbomix}
Y \ \distr \ U \; Z 
\end{equation}
where $U$ is uniform on $[-1,1]$ and $Z \sim (r^{1/r}) W^{1/r}$ with $W \sim$ Ga($1+1/r, 1$), i.e, a gamma distribution with scale 1 and shape parameter $1+1/r$. Hence $Z$ follows a generalized gamma (GG) distribution, specifically $Z \sim $ GG($r, r^{1/r}, 1 + 1/r$). The moments of $Z$ are \citep[see, e.g.,][p.~151]{mom:Kleiber+Kotz:2003}

\[
\mathbb{E}[Z^k] \ = \ \frac{(r ^{1/r})^{1+1/r} \; \Gamma(1 + (k+1)/r)}{\Gamma(1 + 1/r)}, \quad k = 1, 2, \dots .
\]

Now $\lim _{r \to \infty} \mathbb{E}[Z^k] = 1$ for all $k$, and it follows that $Z = r^{1/r} W^{1/r}$ tends to a point mass at 1 by Fr\'{e}chet-Shohat \citep[e.g.,][p.~81]{mom:Galambos:1995}. 
Thus $\lim _{r \to \infty} Y \distr U$, and the density of $\exp(U)$ is given by 

\begin{equation}\label{dtpareto}
f(x) \ = \ \frac{1}{2 \, x}, \quad e^{-1} \leq x \leq e.
\end{equation}

This distribution has compact support, hence it is determined by its moments. \hfill $\square$

\bigskip 
\cite{mom:Lunetta:1963} provides an alternative derivation of the limiting distribution of the generalized normal distribution that analyzes the limit of its characteristic function. However, we prefer the approach involving a mixture representation presented here because it motivates further questions, on which more below. 

\medskip
Interestingly, \cite{mom:Bomsdorf:1977} observed that a distribution of the type described by eq.~(\ref{dtpareto}) occurs as the distribution of prizes in lotteries, hence he calls it the \emph{prize competition distribution}. Among other characteristics he also provides the MGF of this object.

\section{Stieltjes classes for moment-indeterminate \protect\\ generalized lognormal distributions}

The preceding section showed that generalized lognormal distributions of  orders $1 < r < \infty$ are M-indeterminate, by way of an existence proof. To round off the discussion, this section provides explicit examples of distributions that are equivalent, in the sense of having identical moments of all orders, to these indeterminate distributions.

A Stieltjes class -- a term coined by \cite{mom:Stoyanov:2004} -- corresponding to a moment-indeterminate distribution $F$ with density $f$ is a set

\[
\mathcal{S}(f, p) \ = \ \{f_\varepsilon(x) \ | \ f_\varepsilon(x) := f(x)[1 + \varepsilon \ p(x)], \ x \in \hbox{supp}(f), \ \varepsilon \in [-1, 1] \}, 
\]
where $p(x)$ is a perturbation function satisfying $-1 \leq p(x) \leq 1$ and $\mathbb{E}[X^k p(X)] = 0$ for all $k = 0,1,2, \dots$.  

It is possible to obtain Stieltjes classes for the generalized lognormal distributions of orders $1 < r < \infty$ that generalize a recently derived Stieltjes class pertaining to the classical lognormal distribution. The construction of the required Stieltjes classes in the following Theorem is adapted from a construction presented by \citet[][Theorem~3]{mom:Stoyanov+Tolmatz:2005}: 

\begin{theorem}\label{th:scgln}
Suppose $X \sim GLN(\mu, \sigma, r)$ with density $f_r$, 
$(\mu, \sigma, r) \in \mathbb{R} \times \mathbb{R}^+ \times (1, \infty)$. 

\begin{enumerate}
\item The function
\begin{eqnarray}
h_r(x) \ 
= \ \left\{ 
    \begin{array}{ll}
  \sin\{ ( x - 1 )^{1/4} \} 
 \exp \left\{ \frac{1}{r \sigma^r } | \ln x - \mu| ^r  + \ln x 
              - ( x - 1 )^{1/4} 
      \right\}, & x > 1, \label{glnpert} \\
   0, & x \leq 1, 
   \end{array} \right.
\end{eqnarray}
is bounded on $\mathbb{R}^+$ for all $(\mu, \sigma, r) \in \mathbb{R} \times \mathbb{R}^+ \times (1, \infty)$, with $\mathbb{E}[X^k h_r(X)] = 0$ for all $k=0,1,2,\dots$.

\item $p_r := h_r /H_r$, with $H_r:=\sup_x |h_r(x)|$, defines a perturbation corresponding to $f_r$. 

\item The family of functions $f_{r,\varepsilon} (x) = f_r(x) [1 + \varepsilon \ p_r(x)]$, $\varepsilon \in [-1, 1]$, defines a Stieltjes class comprising distributions whose moments are identical to those of $f_r$ for any $\varepsilon \in [-1, 1]$.
\end{enumerate}
\end{theorem}

\noindent
{\sc Proof.} \ The function $h_r$ is continuous on $(1, \infty)$, with $\lim_{x \to 1^+} h_r(x) < \infty $ and $\lim_{x\to \infty} h_r(x) = 0$, hence $h_r$ is bounded on $\mathbb{R}^+$. 

By construction, with $C_r > 0$ the normalizing constant of $f_r$,
\begin{eqnarray*}
\int_0 ^\infty x^k h_{r} (x) f_r(x) \ \hbox{d}x 
&=& C_r \ \int_1 ^\infty x^k \sin\{ ( x - 1 )^{1/4} \} \exp \left\{ - ( x - 1 )^{1/4} \right\} \hbox{d}x \\
&=& C_r \ \int_0 ^\infty (x+1)^k \sin\{ x^{1/4} \} \exp \left\{ - x ^{1/4} \right\} \hbox{d}x \\
&=& C_r \ \sum_{j=0}^k {k \choose j} \int_0 ^\infty x^{k-j} \sin\{ x^{1/4} \} \exp \left\{ - x ^{1/4} \right\} \hbox{d}x
= 0
\end{eqnarray*}
for $k = 0,1,2, \dots$ in view of Lemma 1 of \cite{mom:Stoyanov+Tolmatz:2005} and the fact that 
\[
\int_0 ^\infty x^k \sin\{ x^{1/4} \} \exp \left\{ - x ^{1/4} \right\} \hbox{d}x = 0 , \quad k = 0,1,2, \dots.
\]
This proves (a). 

Since $H_r := \sup_x |h_r(x)| < \infty$  we may set $p_r(x) = h_r(x) / H_r$, assuring $|p_r(x)| \leq 1$ for all $x$. This gives (b). Finally, $f_{r,\varepsilon} (x) = f_r(x) [1 + \varepsilon \ p_r(x)]$ defines a density for any $\varepsilon \in [-1, 1]$, which gives (c). \hfill $\square$

\bigskip
It should be noted that the construction of \cite{mom:Stoyanov+Tolmatz:2005} is somewhat more general, in that the kernel $k(x) := (x-1)^{1/4}$ used here may be generalized to a three-parameter family of kernels defined by $k(x;\xi, \delta, \beta) := (\delta x - \xi)^\beta \tan (\pi \beta)$, where $(\xi, \delta, \beta) \in \mathbb{R}^+ \times \mathbb{R}^+ \times (0,1/2)$. Thus amending the kernel in this manner defines a four-parameter family of perturbations $p_r(x;\xi, \delta, \beta)$ leading to Stieltjes classes that generalize the three-parameter family of Stieltjes classes for the classical lognormal distribution derived by \cite{mom:Stoyanov+Tolmatz:2005}. However, the Stieltjes class presented above already provides infinitely many distributions whose moments coincide with those of the generalized lognormal distribution.

In (\ref{glnpert}), the choice of $\beta = 1/4$ was made because it is related to one of the classical examples of an M-indeterminate distribution that dates back to the pioneering work of \cite{mom:Stieltjes:1894}. Stieltjes considered the case where $\xi = 0$ and the perturbation $h(x) = \sin(x^{1/4})$, $x > 0$, used in the proof of part (a) of Theorem \ref{th:scgln}; it pertains to a certain generalized gamma distribution. Moreover, a shift $\xi > 0$ is needed in (\ref{glnpert}), as otherwise the resulting object would exhibit a singularity at the origin, see also the discussion in \citet[][Section~4]{mom:Stoyanov+Tolmatz:2005}.

\section{Further discussion and concluding remarks}

The paper exhibited a family of distributions, occurring in economics and finance, that generalizes the lognormal distribution, the classical example of a moment-indeterminate distribution. It emerged that a large subfamily consists of moment-indeterminate distributions, but also that not all members share this property of the lognormal, for different reasons: some tails are so heavy that not enough moments exist, while a limiting case corresponds to a light-tailed distribution with compact support.

It may, therefore, be asked to what extent it is possible to characterize the generalized lognormal distributions with $r=1$, i.e. the log-Laplace distributions,%\footnote{This question was raised by an anonymous reviewer.} 
for which $\mathbb{E}[X^k] < \infty$ iff $|k| < 1/\sigma$. 
If one leaves the classical setting of the moment problem characterizations in terms of certain moments are possible. First, Th.~1 of \cite{mom:Lin:1992} implies that characterizations in terms of fractional moments are feasible: for a sequence $\{k_n \; | \; 0 < k_n < 1/\sigma ; n \in \mathbb{N}\}$ of positive and distinct numbers converging to some $k_0 \in (0, 1/\sigma)$, the sequence $\{ \mathbb{E}[X^{k_n}] \; | \; n \in \mathbb{N} \}$ of fractional moments characterizes the distribution. Second, observe that for $r=1$ the first moment exists iff $\sigma < 1$. It is well known that existence of the first moment permits characterization of the underlying distribution in terms of the triangular array of first moments of the associated order statistics, $\{ \mathbb{E}[X_{k:n}] \; | \; k = 1,2, \ldots, n ; n \in \mathbb{N} \}$, where $ X_{1:n} \leq X_{2:n} \leq \ldots \leq X_{n:n}$ are the order statistics in a sample of size $n$. In fact, certain subsets of this array are already sufficient, see \cite{mom:Huang:1989} for a review. Such characterizations are meaningful in applications to income distribution \citep{mom:Kleiber+Kotz:2002}, one of the fields where the generalized lognormal distribution has been employed. Note also that both characterizations, via fractional moments as well as via moments of order statistics, are available for all generalized lognormal distributions with $r > 1$ since moments of arbitrary order exist in that case.

It is natural to ask about M-determinacy of the more widely known distribution of $\ln X$, the generalized error or Subbotin distribution (\ref{generror}). This is a Hamburger moment problem. The answer is already available in the literature, although not in a probabilistic setting: the family of generalized error distributions also admits M-indet examples, namely for $r < 1$, and a Stieltjes class is given in \citet[][p.~22]{mom:Shohat+Tamarkin:1950}.

It is also known that for some M-determinate distributions power transformations lead to M-indeterminacy and vice versa \citep[e.g.][]{mom:Stoyanov:1997}. The standard example is the generalized gamma distribution. For $X \sim$ GLN($\mu, \sigma, r$), it is easily seen that $X^p \sim$ GLN($p\mu, p\sigma, r$) for all $p > 0$, showing that the distribution is closed under power transformations. Hence this well-known property of the classical lognormal distribution extends to the generalized version (\ref{genlnpdf}). Consequently, consideration of power transformations does not lead to new insights regarding the moment problem here.

However, it might be worthwhile to further explore aspects of the mixture representation (\ref{subbomix}). This representation is a special case of a general mixture representation for unimodal distributions known as Khinchine's theorem. The exponentiated version states that $\exp(Y) = \exp(U Z)$, i.e. a random variable following a generalized lognormal distribution can be obtained as the exponential of the product of a uniform and a transformed gamma random variable. It would be interesting to characterize the set of mixing distributions $F_Z$ leading to indeterminate log-unimodal distributions.

%\section*{Acknowledgments}

%I am grateful to Thomas Zehrt for helpful discussions and to an anonymous reviewer for a careful reading of an earlier draft. 

\bibliography{glnmom}

\begin{thebibliography}{29}
\newcommand{\enquote}[1]{``#1''}
\expandafter\ifx\csname natexlab\endcsname\relax\def\natexlab#1{#1}\fi

\bibitem[\protect\citeauthoryear{Alfarano, Milakovi\'{c}, Irle, and
  Kauschke}{Alfarano et~al.}{2012}]{mom:Alfarano+Milakovic+Irle:2012}
\textsc{Alfarano, S., M.~Milakovi\'{c}, A.~Irle, and J.~Kauschke} (2012):
  \enquote{A Statistical Equilibrium Model of Competitive Firms,} \emph{Journal
  of Economic Dynamics and Control}, 36, 136--149.

\bibitem[\protect\citeauthoryear{Bomsdorf}{Bomsdorf}{1977}]{mom:Bomsdorf:1977}
\textsc{Bomsdorf, E.} (1977): \enquote{The Prize-Competition Distribution: A
  Particular {$L$}-Distribution as a Supplement to the {P}areto Distribution,}
  \emph{Statistical Papers}, 18, 254--264.

\bibitem[\protect\citeauthoryear{Box and Tiao}{Box and
  Tiao}{1973}]{mom:Box+Tiao:1973}
\textsc{Box, G. E.~P. and G.~Tiao} (1973): \emph{Bayesian Inference in
  Statistical Analysis}, Reading, MA: Addison-Wesley.

\bibitem[\protect\citeauthoryear{Brunazzo and Pollastri}{Brunazzo and
  Pollastri}{1986}]{mom:Brunazzo+Pollastri:1986}
\textsc{Brunazzo, A. and A.~Pollastri} (1986): \enquote{Proposta di una nuova
  distribuzione: la lognormale generalizzata,} in \emph{Scritti in Onore di
  Francesco Brambilla}, Milano: Bocconi Comunicazioni, vol.~1, pp. 57--83.

\bibitem[\protect\citeauthoryear{Devroye}{Devroye}{1986}]{mom:Devroye:1986}
\textsc{Devroye, L.} (1986): \emph{Non-Uniform Random Number Generation}, New
  York: Springer-Verlag.

\bibitem[\protect\citeauthoryear{Galambos}{Galambos}{1995}]{mom:Galambos:1995}
\textsc{Galambos, J.} (1995): \emph{Advanced Probability Theory}, New York:
  Marcel Dekker, 2nd ed.

\bibitem[\protect\citeauthoryear{Gut}{Gut}{2002}]{mom:Gut:2002}
\textsc{Gut, A.} (2002): \enquote{On the Moment Problem,} \emph{Bernoulli}, 8,
  407--421.

\bibitem[\protect\citeauthoryear{Heyde}{Heyde}{1963}]{mom:Heyde:1963}
\textsc{Heyde, C.~C.} (1963): \enquote{On a Property of the Lognormal
  Distribution,} \emph{Journal of the Royal Statistical Society, Series B}, 25,
  392--393.

\bibitem[\protect\citeauthoryear{Huang}{Huang}{1989}]{mom:Huang:1989}
\textsc{Huang, J.~S.} (1989): \enquote{Moment Problem of Order Statistics: A
  Review,} \emph{International Statistical Review}, 57, 59--66.

\bibitem[\protect\citeauthoryear{Kleiber}{Kleiber}{2000}]{mom:Kleiber:2000}
\textsc{Kleiber, C.} (2000): \enquote{A Simple Distribution Without Any
  Moments,} \emph{The Mathematical Scientist}, 25, 59--60.

\bibitem[\protect\citeauthoryear{Kleiber and Kotz}{Kleiber and
  Kotz}{2002}]{mom:Kleiber+Kotz:2002}
\textsc{Kleiber, C. and S.~Kotz} (2002): \enquote{A Characterization of Income
  Distributions in Terms of Generalized {G}ini Coefficients,} \emph{Social
  Choice and Welfare}, 19, 789--794.

\bibitem[\protect\citeauthoryear{Kleiber and Kotz}{Kleiber and
  Kotz}{2003}]{mom:Kleiber+Kotz:2003}
---\hspace{-.1pt}---\hspace{-.1pt}--- (2003): \emph{Statistical Size
  Distributions in Economics and Actuarial Sciences}, Hoboken, NJ: John Wiley
  \& Sons.

\bibitem[\protect\citeauthoryear{Kotz}{Kotz}{1975}]{mom:Kotz:1975}
\textsc{Kotz, S.} (1975): \enquote{Multivariate Distributions at a Cross-Road,}
  in \emph{Statistical Distributions in Scientific Work}, ed. by G.~P. Patil,
  S.~Kotz, and J.~K. Ord, Dordrecht: D. Reidel Publishing Company, vol.~1,
  247--270.

\bibitem[\protect\citeauthoryear{Lin}{Lin}{1992}]{mom:Lin:1992}
\textsc{Lin, G.~D.} (1992): \enquote{Characterizations of Distributions via
  Moments,} \emph{Sankhy\={a}}, A 54, 128--132.

\bibitem[\protect\citeauthoryear{Lin and Stoyanov}{Lin and
  Stoyanov}{2009}]{mom:Lin+Stoyanov:2009}
\textsc{Lin, G.~D. and J.~Stoyanov} (2009): \enquote{The Logarithmic
  Skew-Normal Distributions are Moment-Indeterminate,} \emph{Journal of Applied
  Probability}, 46, 909--916.

\bibitem[\protect\citeauthoryear{Lunetta}{Lunetta}{1963}]{mom:Lunetta:1963}
\textsc{Lunetta, G.} (1963): \enquote{Di una generalizzazione dello schema
  della curva normale,} \emph{Annali della Facolt\`{a} di Economia e Commercio
  di Palermo}, 17, 237--244.

\bibitem[\protect\citeauthoryear{Nelson}{Nelson}{1991}]{mom:Nelson:1991}
\textsc{Nelson, D.~B.} (1991): \enquote{Conditional Heteroskedasticity in Asset
  Returns: {A} New Approach,} \emph{Econometrica}, 59, 347--370.

\bibitem[\protect\citeauthoryear{Pakes, Hung, and Wu}{Pakes
  et~al.}{2001}]{mom:Pakes+Hung+Wu:2001}
\textsc{Pakes, A.~G., W.-L. Hung, and J.-W. Wu} (2001): \enquote{Criteria for
  the Unique Determination of Probability Distributions by Moments,}
  \emph{Australian \& New Zealand Journal of Statistics}, 43, 101--111.

\bibitem[\protect\citeauthoryear{Shohat and Tamarkin}{Shohat and
  Tamarkin}{1950}]{mom:Shohat+Tamarkin:1950}
\textsc{Shohat, J.~A. and J.~D. Tamarkin} (1950): \emph{The Problem of
  Moments}, Providence, RI: American Mathematical Society, revised ed.

\bibitem[\protect\citeauthoryear{Stieltjes}{Stieltjes}{1894/1895}]{mom:Stieltjes:1894}
\textsc{Stieltjes, T.~J.} (1894/1895): \enquote{Recherches sur les fractions
  continues,} \emph{Annales de la Facult\'{e} des Sciences de Toulouse}, 8/9,
  1--122, 1--47.

\bibitem[\protect\citeauthoryear{Stoyanov}{Stoyanov}{1997}]{mom:Stoyanov:1997}
\textsc{Stoyanov, J.} (1997): \emph{Counterexamples in Probability},
  Chichester: John Wiley \& Sons, 2nd ed.

\bibitem[\protect\citeauthoryear{Stoyanov}{Stoyanov}{2000}]{mom:Stoyanov:2000}
---\hspace{-.1pt}---\hspace{-.1pt}--- (2000): \enquote{Krein Condition in
  Probabilistic Moment Problems,} \emph{Bernoulli}, 6, 939--949.

\bibitem[\protect\citeauthoryear{Stoyanov}{Stoyanov}{2004}]{mom:Stoyanov:2004}
---\hspace{-.1pt}---\hspace{-.1pt}--- (2004): \enquote{Stieltjes Classes for
  Moment-Indeterminate Probability Distributions,} \emph{Journal of Applied
  Probability}, 41A, 281--294.

\bibitem[\protect\citeauthoryear{Stoyanov and Tolmatz}{Stoyanov and
  Tolmatz}{2005}]{mom:Stoyanov+Tolmatz:2005}
\textsc{Stoyanov, J. and L.~Tolmatz} (2005): \enquote{Methods for Constructing
  {S}tieltjes Classes for {M}-Indeterminate Probability Distributions,}
  \emph{Applied Mathematics and Computation}, 165, 669--685.

\bibitem[\protect\citeauthoryear{Subbotin}{Subbotin}{1923}]{mom:Subbotin:1923}
\textsc{Subbotin, M.~T.} (1923): \enquote{On the Law of Frequency of Error,}
  \emph{Mathematicheskii Sbornik}, 31, 296--301.

\bibitem[\protect\citeauthoryear{Taylor}{Taylor}{2005}]{mom:Taylor:2005}
\textsc{Taylor, S.~J.} (2005): \emph{Asset Price Dynamics, Volatility, and
  Prediction}, Princeton, NJ: Princeton University Press.

\bibitem[\protect\citeauthoryear{Vianelli}{Vianelli}{1963}]{mom:Vianelli:1963}
\textsc{Vianelli, S.} (1963): \enquote{La misura della variabilit\`{a}
  condizionata in uno schema generale delle curve normali di frequenza,}
  \emph{Statistica}, 23, 447--474.

\bibitem[\protect\citeauthoryear{Vianelli}{Vianelli}{1982}]{mom:Vianelli:1982a}
---\hspace{-.1pt}---\hspace{-.1pt}--- (1982): \enquote{Sulle curve lognormali
  di ordine $r$ quali famiglie di distribuzioni di errori di proporzione,}
  \emph{Statistica}, 42, 155--176.

\bibitem[\protect\citeauthoryear{Vianelli}{Vianelli}{1983}]{mom:Vianelli:1983}
---\hspace{-.1pt}---\hspace{-.1pt}--- (1983): \enquote{The Family of Normal and
  Lognormal Distributions of Order $r$,} \emph{Metron}, 41, 3--10.

\end{thebibliography}

\end{document}